\documentclass{amsart}
\usepackage{amssymb}
\usepackage{upref}

\usepackage[all,cmtip]{xy}

\usepackage[lite,initials,shortalphabetic]{amsrefs}

\providecommand{\texorpdfstring}[2]{#1}

\makeatletter
  \renewcommand{\p@enumi}{\thesubsection}
\makeatother

\newenvironment{resumeenumerate}[1]
{\begin{enumerate}
 \setcounter{enumi}{#1}
 \addtocounter{enumi}{-1}
}
{\end{enumerate}
}

\newenvironment{lettered}
{\begin{list}{\thelettercounter.}
 {\usecounter{lettercounter}\def\makelabel##1{\hss\llap{##1}}}
}
{\end{list}
}
\newcounter{lettercounter}
\renewcommand{\thelettercounter}{\Alph{lettercounter}}

\theoremstyle{plain}

\makeatletter
\newcommand{\emsection}[1]{%
  \par
  \addpenalty\@secpenalty
  \vskip 6 pt plus 9 pt
  \emph{#1.}\nobreak\enspace\ignorespaces
}
\makeatother

\newcommand{\intro}{%
  \goodbreak
  \vskip 6 pt plus 9 pt
}

\numberwithin{equation}{subsection}

\newcommand{\Period}{\rlap{\enspace .}}

\newcommand{\cat}[1]{\boldsymbol{#1}}

\newcommand{\bs}{\boldsymbol}

\newcommand{\RelCat}{\mathbf{RelCat}}

\newcommand{\SCat}{\mathbf{SCat}}

\newcommand{\SOCat}{\mathbf{S}O\text{-}\mathbf{Cat}}

\newcommand{\RelSCat}{\mathbf{RelSCat}}

\newcommand{\DK}{\mathbf{DK}}

\DeclareMathOperator{\Ho}{Ho}
\DeclareMathOperator{\colim}{colim}

\DeclareMathOperator{\Rel}{Rel}
\DeclareMathOperator{\Loc}{L}

\newcommand{\uH}{^{\mathrm{H}}}
\newcommand{\uc}{^{\mathrm{c}}}
\newcommand{\ucp}{^{\mathrm{c}'}}
\newcommand{\Fstar}{\mathrm{F}_{*}}

\newcommand{\id}{\mathbf{id}}

\newcommand{\iso}{\approx}

\newcommand{\pushout}[3]{#1\mathbin{\mathord{\smallcoprod}_{#2}}#3}

\newcommand{\smallcoprod}{\mathchoice{\mathbin\amalg}%
               {\mathbin\amalg}%
               {{\scriptscriptstyle\mathbin{\amalg}}}%
               {{\scriptscriptstyle\mathbin{\amalg}}}}

\newcommand{\union}{\cup}

\begin{document}

\title{A characterization of simplicial localization functors}

\author{C. Barwick}
\address{Department of Mathematics, Massachusetts Institute of
  Technology, Cambridge, MA 02139}
\email{clarkbar@gmail.com}

\author{D.M. Kan}
\address{Department of Mathematics, Massachusetts Institute of
  Technology, Cambridge, MA 02139}

\date{\today}

\begin{abstract}
  We characterize simplicial localization functors among relative functors from relative categories to simplicial categories as any choice of homotopy inverse to the delocalization functor of Dwyer and the second author.
\end{abstract}

\maketitle

\section{An overview}
\label{sec:Ovrvw}

We start with some preliminaries.

\subsection{Relative categories}
\label{sec:RelCat}

As in \cite{BK} we denote by $\RelCat$ the category of (small)
relative categories and relative functors between them, where by a
\textbf{relative category} we mean a pair $(\cat C, \cat W)$
consisting of a category $\cat C$ and a subcategory $\cat W \subset
\cat C$ which contains all the objects of $\cat C$ and their identity
maps and of which the maps will be referred to as \textbf{weak
  equivalences} and where by a \textbf{relative functor} between two
such relative categories we mean a weak equivalence preserving
functor.

\subsection{Homotopy equivalences between relative categories}
\label{sec:HomEqRlCt}

A relative functor $f\colon \cat X \to \cat Y$ between two relative
categories \eqref{sec:RelCat} is called a \textbf{homotopy
  equivalence} if there exists a relative functor $g\colon \cat Y \to
\cat X$ (called a \textbf{homotopy inverse} of $f$) such that the
compositions $gf$ and $fg$ are naturally weakly equivalent (i.e.\ can
be connected by a finite zigzag of natural weak equivalences) to the
identity functors of $\cat X$ and $\cat Y$ respectively.

\subsection{DK-equivalences}
\label{sec:DKeq}

A map in the category $\SCat$ of \emph{simplicial categories} (i.e.\
categories enriched over simplicial sets) is \cite{Be} called a
DK-\textbf{equivalence} if it induces \emph{weak equivalences} between
the simplicial sets involved and an \emph{equivalence of categories}
between their \emph{homotopy categories}, i.e.\ the categories
obtained from them by replacing each simplicial set by the set of its
components.

Furthermore a map in $\RelCat$ will similarly be called a
DK\textbf{equivalence} if its image in $\SCat$ is so under the
\emph{hammock localization functor} \cite{DK2}
\begin{displaymath}
  L^{H}\colon \RelCat\longrightarrow \SCat
\end{displaymath}
(or of course the naturally DK-equivalent functors $\RelCat \to \SCat$
considered in \cite{DK1} and \cite{DHKS}*{35.6}).

We will denote by both
\begin{displaymath}
  \DK\subset \SCat
  \qquad\text{and}\qquad
  \DK\subset \RelCat
\end{displaymath}
the subcategories consisting of these DK-equivalences.

\intro
Next we define what we mean by
\subsection{Simplicial localization functors}
\label{sec:SpLcFn}

In defining DK-equivalences in $\RelCat$ \eqref{sec:DKeq} we used the
hammock localization functor and not one of the other DK-equivalent
functors mentioned because, for our purposes here it seemed to be the
more convenient one.  However in other situations the others are more
convenient and it therefore makes sense to define in general a
\textbf{simplicial localization functor} as any functor $\RelCat \to
\SCat$ which is naturally DK-equivalent to the functors mentioned
above \eqref{sec:DKeq}.

\intro
We also need
\subsection{The relativization functor}
\label{sec:RlFunc}

In contrast with the situation mentioned in \ref{sec:SpLcFn} there
\emph{is} a preferred choice for a \textbf{relativization functor}
\begin{displaymath}
  R\colon \SCat \longrightarrow \RelCat
\end{displaymath}
which is a kind of inverse of the simplicial localization functor,
namely the \emph{delocalization} mentioned in \cite{DK3}*{2.5} which
assigns to an object $\cat A \in \SCat$ its relative \emph{flattening}
which is the relative category which
consists of
\begin{enumerate}
\item a category which is the Grothendieck construction on $\cat A$,
  where $\cat A$ is considered as a simplicial diagram of categories,
  and
\item its subcategory obtained by applying the same construction to
  the subobject of $\cat A$ which consists of its objects only.
\end{enumerate}

\intro
Our main result then is
\subsection{Theorem}
\label{thm:CharSimpLoc}
\begin{em}
  A relative functor
  \begin{displaymath}
    (\RelCat, \DK) \longrightarrow (\SCat, \DK)
  \end{displaymath}
  is a simplicial localization functor \eqref{sec:RlFunc} iff it is a
  homotopy inverse \eqref{sec:HomEqRlCt} of the realization functor
  \eqref{sec:RlFunc}
  \begin{displaymath}
    \Rel\colon (\SCat, \DK) \longrightarrow (\RelCat, \DK) \Period
  \end{displaymath}
\end{em}

\intro
We end with some
\subsection{Comments on the proof of  \ref{thm:CharSimpLoc}}
\label{sec:CmntPrf}

The proof of theorem~\ref{thm:CharSimpLoc} heavily involves some of
the results of \cite{DK1} and \cite{DK3}*{2.5} and we therefore first
(in \S \ref{sec:Prelim}) review some of the results of these papers.

In \S \ref{sec:PrfThm} we then actually prove theorem
\ref{thm:CharSimpLoc}.  It turns out however that in addition to the
results mentioned in \S \ref{sec:Prelim} we need a property of the
hammock localization of which we will give two proofs.  The first is a
very short one based on a remark of Toen and Vezzosi \cite{TV}*{2.2.1}
involving the \emph{homotopy category} of $\SCat$.  The other, which
is due to Bill Dwyer, relies heavily on \cite{DK1} and \cite{DK2} and
is longer, but has the ``advantage'' of taking place in the
\emph{model category} itself.

\section{Preliminaries}
\label{sec:Prelim}

In preparation for the proof (in \S \ref{sec:PrfThm}) of theorem
\ref{thm:CharSimpLoc} we review here some of the results of
\cite{DK1}, \cite{DK2} and \cite{DK3}*{2.3} which will be needed.

\subsection{The hammock localization}
\label{sec:HmkLoc}

In the proof of \ref{thm:CharSimpLoc} we will make extensive use of
the hammock localization $\Loc\uH$ of \cite{DK2} because, unlike the
other simplicial localization functors, it has the property that every
relative category $(\cat C, \cat W)$ comes with a \emph{natural
  embedding} $\cat C \to \Loc\uH(\cat C, \cat W)$.

\subsection{The category \texorpdfstring{$\RelSCat$}{RelSCat}}
\label{sec:RelSCat}

This will be the category which has as its objects the pairs $(\cat A,
\cat U)$ where $\cat A \in \SCat$ and $\cat U \subset \cat A$ is a
subobject which contains all the objects of $\cat A$.

One then can consider $\RelCat$ as a full subcategory of $\RelSCat$
and \cite{DK2}*{2.5} extend the functor $\Loc\uH\colon \RelCat \to
\SCat$ to a functor $\Loc\uH\colon \RelSCat \to \SCat$ by sending an
object of $\RelSCat$ to the diagonal of the bisimplicial set obtained
from it by dimensionwise application of the hammock localization.

\intro
To deal with \cite{DK1} and \cite{DK3}*{2.5} it will be convenient to
introduce a notion of
\subsection{Neglectable categories}
\label{sec:Neglect}

Given an object $(\cat A, \cat U) \in \RelSCat$ \eqref{sec:RelSCat} we
will say that $\cat U$ is \emph{neglectable} in $\cat A$ if every map
of $\cat U$ goes to an \emph{isomorphism} in $\pi_{0}\cat A$.

\subsection{Some results from \texorpdfstring{\cite{DK1}}{}}
\label{sec:SmRslts}

\cite{DK1}*{3.4 and 5.1} imply
\begin{enumerate}
\item
  \begin{em}
    Let $\cat A$ be a category, let $\cat U$ and $\cat V\subset \cat
    A$ be subcategories which contain all the objects of $\cat A$ and
    let $\cat U\union \cat V \subset \cat A$ denote the subcategory
    spanned by $\cat U$ and $\cat V$ and assume that $\cat V$ is
    neglectable in $\Loc\uH(\cat A, \cat U)$.  Then the induced map
    \begin{displaymath}
      \Loc\uH(\cat A, \cat U) \longrightarrow
      \Loc\uH(\cat A, \cat U\union \cat V) \in \SCat
    \end{displaymath}
    is a DK-equivalence.
  \end{em}
\end{enumerate}

Similarly \cite{DK1}*{6.4} implies
\begin{resumeenumerate}{2}
\item
  \begin{em}
    Let $(\cat B, \cat V) \in \RelSCat$ be such that $\cat V$ is
    neglectable in $\cat B$.  Then the induced map (\ref{sec:HmkLoc}
    and \ref{sec:RelSCat})
    \begin{displaymath}
      \cat B \longrightarrow \Loc\uH(\cat B,\cat V) \in \SCat
    \end{displaymath}
    is a DK-equivalence.
  \end{em}
\end{resumeenumerate} 

\intro
We end with a brief review of
\subsection{The relativization functor \texorpdfstring{\cite{DK3}*{2.5}}{}}
\label{sec:ReltvFunc}

The \textbf{relativization functor} is the functor
\begin{displaymath}
  \Rel\colon \SCat \longrightarrow \RelCat
\end{displaymath}
which sends an object $\cat A \in \SCat$ to the object $(b\cat A,
b\id) \in \RelCat$, where $b\cat A$ is the \emph{flattening} of $\cat
A$, i.e.\ the category which has as objects the pairs $(A,n)$, where
$A$ is an object of $\cat A$ and $n$ is an integer $\ge 0$ and which
has as maps $(A_{1}, n_{1}) \to (A_{2}, n_{2})$ the pairs $(a,q)$
where $a$ is a map $A_{1}\to A_{2} \in \cat A_{n_{2}}$ and $q$ is a
simplicial operator from dimension $n_{1}$ to dimension $n_{2}$ and
$b\id\in b\cat A$ is the subcategory consisting of the maps $(a,q)$
for which $a$ is an identity map.

It then was noted in \cite{DK3}*{2.5} that, for every object $\cat A
\in \SCat$, there exists an object $\overline{\cat A} \in \SCat$ with
the same object set as $b\cat A$ with the following properties:
\begin{em}
  \begin{enumerate}
  \item There is a natural monomorphism $\cat A \to \overline{\cat A}$
    which is a DK-equivalence.
  \item There is a natural embedding $b\cat A \to \overline{\cat A}$
    with the property that (if the image of $b\id \in b\cat A$ in
    $\overline{\cat A}$ is also denoted by $b\id$) the induced map
   \begin{displaymath}
     \Loc\uH(b\cat A, b\id) \longrightarrow
     \Loc\uH(\overline{\cat A}, b\id) \in \SCat
   \end{displaymath}
   is a DK-equivalence.
 \item $b\id$ is neglectable in $\overline{\cat A}$
   \eqref{sec:Neglect} and hence the embedding (\ref{sec:HmkLoc} and
   \ref{sec:RelSCat})
   \begin{displaymath}
     \overline{\cat A}\longrightarrow
     \Loc\uH(\overline{\cat A}, b\id) \in \SCat
   \end{displaymath}
   is a DK-equivalence.
  \end{enumerate}
\end{em}
It follows that
\begin{em}
  \begin{resumeenumerate}{4}
  \item $\cat A$ and $\Loc\uH\Rel\cat A$ can be connected by the
    natural zigzag of DK-equivalences
    \begin{displaymath}
      \cat A \longrightarrow
      \overline{\cat A} \longrightarrow
      \Loc\uH(\overline{\cat A}, b\id) \longleftarrow
      \Loc\uH(b\cat A, b\id) =
      \Loc\uH\Rel\cat A
    \end{displaymath}
  \end{resumeenumerate}
\end{em}
which in turn implies that
\begin{em}
  \begin{resumeenumerate}{5}
  \item  $\Rel$ is a relative functor \eqref{sec:DKeq}
    \begin{displaymath}
      \Rel\colon (\SCat, \DK) \longrightarrow (\RelCat, \DK) \Period
    \end{displaymath}
  \end{resumeenumerate}
\end{em}

\section{A proof of theorem \ref{thm:CharSimpLoc}}
\label{sec:PrfThm}

To prove theorem \ref{thm:CharSimpLoc} is suffices, in view of
\ref{sec:ReltvFunc}(iv) and (v), to prove
\subsection{Proposition}
\label{prop:DKeq}
\begin{em}
  Every object $(\cat C, \cat W) \in \RelCat$ is naturally
  DK-equivalent to $\Rel\Loc\uH(\cat C, \cat W)$.
\end{em}

\emsection{Proof}
Consider the commutative diagram in $\RelSCat$ \eqref{sec:RelSCat}
\begin{displaymath}
  \xymatrix@C=0em{
    {(\cat C,\cat W)} \ar[dr]^{a} \ar[dd]_{c}
    && {\bigl(bL^{H}(\cat C,\cat W),b\bs{\id}\bigr)
      = RL^{H}(\cat C,\cat W)} \ar[dd]^{e} \ar[dl]_{b}\\
    & {\bigl(bL^{H}(\cat C,\cat W), b\id\union \cat W\bigr)}
    \ar[dd]^{d}\\
    {\bigl(L^{H}(\cat C,\cat W),\cat W\bigr)} \ar[dr]_{f}
    && {\bigl(\overline{L^{H}(\cat C,\cat W)}, b\id\bigr)}
    \ar[dl]^{g}\\
    & {\bigl(\overline{L^{H}(\cat C,\cat W)}, b\id\union \cat W\bigr)}
  }
\end{displaymath}
in which
\begin{itemize}
\item $c$ is as in \ref{sec:HmkLoc}, $f$ is as in
  \ref{sec:ReltvFunc}(i), $d$ and $e$ are as in
  \ref{sec:ReltvFunc}(ii) and $a$ is the unique map such that $da =
  fc$, and
\item the symbol $\union$ is as in \ref{sec:SmRslts}(i) and, in the
  formulas which involve two $\cat W$'s, the second $\cat W$ is the
  image of the $\cat W$ in the upper left $(\cat C, \cat W)$.
\end{itemize}

Then it suffices to show that $a$ and $b$ are DK-equivalences in
$\RelCat$ or equivalently that $\Loc\uH a$ and $\Loc\uH b$ are
DK-equivalences in $\SCat$.

This is done as follows:

The map $f$ admits a factorization
\begin{displaymath}
  \bigl(\Loc\uH(\cat C, \cat W), \cat W\bigr)
  \xrightarrow{\enspace x\enspace}
  \bigl(\overline{\Loc\uH(\cat C, \cat W)}, \cat W\bigr)
  \xrightarrow{\enspace y\enspace}
  \bigl(\overline{\Loc\uH(\cat C, \cat W)}, b\id\union \cat W\bigr)
\end{displaymath}
in which clearly $\cat W$ is neglectable \eqref{sec:Neglect} in
$\Loc\uH(\cat C, \cat W)$ and hence, in view of \ref{sec:SmRslts}(ii)
and \ref{sec:ReltvFunc}(i), $\Loc\uH x$ is a DK-equivalence and $\cat
W$ is neglectable in $\overline{\Loc\uH(\cat C, \cat W)}$.  It follows
that (\ref{sec:ReltvFunc}(iii)) $\Loc\uH y$ is a DK-equivalence and
hence (\ref{sec:SmRslts}(ii)) so is $\Loc\uH g$.

Furthermore, in view of \ref{sec:ReltvFunc}(ii), $\Loc\uH e$ is a
DK-equivalence and consequently $\Loc\uH d$ is a DK-equivalence and
$\cat W$ is neglectable in $\Loc\uH\bigl(b\Loc\uH(\cat C, \cat W),
b\id\bigr)$ which implies that $\Loc\uH b$ is a DK-equivalence.

It thus remains to prove that $\Loc\uH a$ is a DK-equivalence, but
this now follows from
\subsection{Proposition}
\label{prop:LocDKeq}
\begin{em}
  $\Loc\uH c\colon \Loc\uH(\cat C, \cat W) \to 
  \Loc\uH\bigl(\Loc\uH(\cat C, \cat W), \cat W\bigr)$ is a
  DK-equivalence.
\end{em}

We will give two proofs of this proposition.

The first is short and is based on a remark by Toen and Vezzosi
involving the \emph{homotopy category} $\Ho\SCat$ of $\SCat$.

The other, due to Bill Dwyer, is longer but takes place inside the
\emph{model category} $\SOCat$ of the simplicial categories with a
fixed object set $O$ (in this case the object set of $\cat C$).

They both involve the commutative diagram
\begin{enumerate}
\item $\vcenter{
    \xymatrix{
      {\cat C} \ar[r] \ar[d]
      & {\Loc\uH(\cat C, \cat W)} \ar[d]\\
      {\Loc\uH(\cat C, \cat W)} \ar[r]^-{\Loc c}
      & {\Loc\uH\bigl(\Loc\uH(\cat C, \cat W), \cat W\bigr)}
    }
  }$\\
  in which the unmarked maps are as in \ref{sec:HmkLoc}, and
\item in which, in view of \ref{sec:SmRslts}(ii) the right hand map is
  a DK-equivalence.
\end{enumerate}

\subsection{The short proof}
\label{sec:ShrtPrf}

In view of \cite{TV}*{2.2.1} the maps at the right and the bottom in
\ref{prop:LocDKeq}(i) have the same image in $\Ho\SCat$ and as
(\ref{prop:LocDKeq}(ii)) the one on the right is a DK-equivalence, so
is the one at the bottom.

\subsection{The longer proof}
\label{sec:LngPrf}

We start with a brief discussion of
\begin{enumerate}
\item \textbf{Homotopy pushouts}

  Given a model category together with a choice of \emph{cofibrant
    replacement functor} and a choice of \emph{functorial
    factorization} of maps into a cofibration followed by a trivial
  fibration, associate with every zigzag $Y \gets X \to Z$ a
  commutative diagram
  \begin{displaymath}
    \xymatrix{
      {Y\uc} \ar[d]_{\sim}
      & {X\uc} \ar[l] \ar[r] \ar[dl] \ar[dr] \ar[dd]_(.75){\sim}
      & {Z\uc} \ar[d]^{\sim}\\
      {Y\ucp} \ar[d]_{\sim}
      && {Z\ucp} \ar[d]^{\sim}\\
      {Y}
      & {X} \ar[l] \ar[r]
      & {Z}
    }
  \end{displaymath}
  as follows.  The pentagon is obtained by applying the cofibrant
  approximation functor and the two triangles by means of the
  functorial factorization.  Consequently the maps indicated $\sim$
  are weak equivalences.

  Then the pushout $\pushout{Y\uc}{X\uc}{Z\uc}$ is a homotopy pushout
  of the zigzag $Y \gets X \to Z$.

  Clearly this construction is functorial in the sense that every
  diagram of the form
  \begin{displaymath}
    \xymatrix{
      {Y_{0}} \ar[d]_{y}
      & {X_{0}} \ar[l] \ar[r] \ar[d]_{x}
      & {Z_{0}} \ar[d]_{z}\\
      {Y_{1}}
      & {X_{1}} \ar[l] \ar[r]
      & {Z_{1}}
    }
  \end{displaymath}
  induces a map
  \begin{displaymath}
    \pushout{Y\uc_{0}}{X\uc_{0}}{Z\uc_{0}} \longrightarrow
    \pushout{Y\uc_{1}}{X\uc_{1}}{Z\uc_{1}}
  \end{displaymath}
  which is a weak equivalence whenever $y$, $x$ and $z$ are.
\end{enumerate}

Next we discuss

\begin{resumeenumerate}{2}
\item \textbf{The original simplicial localization functor $\Loc$
    \cite{DK1}}

  We will work in the model category $\SOCat$ \cite{DK1} of the
  simplicial categories with a fixed object set $O$ (which will be the
  object set of $\cat C$).  The weak equivalences in the model
  structure are the DK-equivalences.

  As \cite{DK1}*{4.1} $\Loc(\cat C, \cat W)$ is the pushout of the zigzag
  \begin{displaymath}
    \Fstar \cat C \longleftarrow \Fstar \cat W \longrightarrow
    \Fstar \cat W[\Fstar \cat W^{-1}]
  \end{displaymath}
  and the map $\Fstar \cat W \to \Fstar \cat C$ is a cofibration, it
  follows from \cite{DK1}*{8.1} that this pushout is also a homotopy
  pushout.  Consequently
  \begin{itemize}
  \item [($*$)] $\Loc(\cat C, \cat W)$ is naturally DK-equivalent to
    \begin{displaymath}
      \pushout{(\Fstar\cat C)\uc}{(\Fstar\cat W)\uc}
      {\bigl(\Fstar \cat W[\Fstar\cat W^{-1}]\bigr)\uc} \Period
    \end{displaymath}
  \end{itemize}
\end{resumeenumerate}

Now we turn to
\begin{resumeenumerate}{3}
\item \textbf{The hammock localization $\Loc\uH$ \cite{DK2}}

  In view of \cite{DK2}*{2.5} the functors $\Loc$ and $\Loc\uH$ are
  naturally DK-equivalent and there exists a diagram of the form
  \begin{displaymath}
    \xymatrix{
      {\Fstar\cat C}
      & {\Fstar\cat W} \ar[l] \ar[r]
      & {\Fstar\cat W[\Fstar\cat W^{-1}]} \ar@{}[r]|-{=}
      & {\Loc(\cat W, \cat W)}\\
      {\Fstar\cat C} \ar[u]^{\iso} \ar[d]
      & {\Fstar\cat W} \ar[l] \ar[u]^{\iso} \ar[r] \ar[d]
      & {\Loc\uH(\Fstar\cat W, \Fstar\cat W)} \ar[u] \ar[d]\\
      {\cat C}
      & {\cat W} \ar[l] \ar[r]
      & {\Loc\uH(\cat W, \cat W)}
    }
  \end{displaymath}
  in which the vertical maps are DK-equivalences.  Hence
  \begin{itemize}
  \item [($*$)] $\Loc\uH(\cat C, \cat W)$ is naturally DK-equivalent
    to
    \begin{displaymath}
      \pushout{\cat C\uc}{\cat W\uc}{\Loc\uH(\cat W,\cat W)\uc}
    \end{displaymath}
    and $\Loc\uH\bigl(\Loc\uH(\cat C,\cat W),\cat W\bigr)$ is
    naturally DK-equivalent to
    \begin{displaymath}
      \cat C\uc \mathbin{\mathord{\smallcoprod}_{\cat W\uc}}
      \Loc\uH(\cat W,\cat W)\uc
      \mathbin{\mathord{\smallcoprod}_{\cat W\uc}}
      \Loc\uH(\cat W,\cat W)\uc
    \end{displaymath}
    in which the two middle maps $\cat W\uc \to \Loc\uH(\cat W,\cat
    W)\uc$ are the same and which therefore is the same as
    \begin{displaymath}
      \cat Q = \colim\left(\vcenter{
          \xymatrix@R=0.5ex{&& {\Loc\uH(\cat W,\cat W)\uc}\\
            {\cat C\uc}
            & {\cat W\uc} \ar[l] \ar[ur] \ar[dr]\\
            && {\Loc\uH(\cat W,\cat W)\uc}}
        }\right)\Period
    \end{displaymath}
  \end{itemize}
  It follows that diagram \ref{prop:LocDKeq}(i) is DK-equivalent to
  the commutative diagram
  \begin{displaymath}
    \xymatrix{
      {\cat C\uc} \ar[r] \ar[d]
      & {\pushout{\cat C\uc}{\cat W\uc}{\Loc\uH(\cat W,\cat W)\uc}}
      \ar[d]_{u}\\
      {\pushout{\cat C\uc}{\cat W\uc}{\Loc\uH(\cat W,\cat W)\uc}}
      \ar[r]^-{v}
      & {\cat Q}
    }
  \end{displaymath}
  in which $u$ is obtained by mapping the zigzag $\cat C\uc \gets \cat
  W\uc \to \Loc\uH(\cat W,\cat W)\uc$ to the upper zigzag in the
  diagram whose colimit is $\cat Q$ and $v$ is obtained by mapping it
  to the lower zigzag.

  In view of \ref{prop:LocDKeq}(ii) the map $u$ is a DK-equivalence
  and as $v = Tu$ where $T\colon \cat Q \to \cat Q$ denotes the
  automorphism which switches the two copies of $\Loc\uH(\cat W,\cat
  W)\uc$, so is the map $v$ and therefore also the desired map
  \begin{displaymath}
    \Loc c\colon \Loc\uH(\cat C,\cat W) \longrightarrow
    \Loc\uH\bigl(\Loc\uH(\cat C,\cat W), \cat W\bigr) \Period
  \end{displaymath}
\end{resumeenumerate}


\end{document}